\documentclass[12pt]{amsart}
\usepackage{amsmath,amssymb,amscd,amsxtra}
\usepackage{latexsym}
\usepackage[dvips]{graphics,epsfig}

\newcommand{\N}{\mathbb{N}}
\newcommand{\R}{\mathbb{R}}
\newcommand{\C}{\mathbb{C}}
\newcommand{\Z}{\mathbb{Z}}
\newcommand{\M}[1]{\mathcal{M}^{#1}}
\renewcommand{\S}{\mathcal{S}}
\newcommand{\di}{\mathcal{S}^{'}}

\newcommand{\nm}[2]{\|{#1}\|_{#2}}

\newcommand{\bigparen}[1]{\bigl(#1\bigr)}

\newcommand{\biggparen}[1]{\biggl(#1\biggr)}
\newcommand{\ip}[2]{\langle#1,#2\rangle}
\newcommand{\bignm}[2]{\bigg{\|}#1\bigg{\|}_{#2}}

\newcommand{\rd}{{\R}^{d}}

\newcommand{\wam}{W(\cF L^1, \ell ^\infty)}
\newcommand{\wamp}{W(\cF L^p, \ell ^\infty)}

\newcommand{\cF}{\mathcal{F}}
\newcommand{\cS}{\mathcal{S}}

\newtheorem{theorem}{Theorem}
\newtheorem{lemma}{Lemma}

\newtheorem{corollary}{Corollary}
\theoremstyle{remark}

\theoremstyle{definition}

\begin{document}

\title{Local well-posedness of nonlinear dispersive equations on
modulation spaces}

\author{\'Arp\'ad B\'enyi}

\address{\'Arp\'ad B\'enyi\\
Department of Mathematics\\
516 High Street\\
Western Washington University\\
Bellingham, WA 98225, USA }

\email{arpad.benyi@wwu.edu}

\author{Kasso A.~Okoudjou}

\address{Kasso A.~Okoudjou\\
Department of Mathematics\\
University of Maryland\\
College Park, MD 20742, USA}

\email{kasso@math.umd.edu}


\date{\today}


\begin{abstract}
By using tools of time-frequency analysis, we obtain some improved
local well-posedness results for the NLS, NLW and NLKG equations
with Cauchy data in modulation spaces $\M{p, 1}_{0,s}$.
\end{abstract}

\maketitle
\pagestyle{myheadings}
\thispagestyle{plain}
\markboth{\'A. B\'enyi and  K. A. Okoudjou}
{NONLINEAR DISPERSIVE EQUATIONS ON MODULATION SPACES}

\section{Introduction and statement of results}

The theory of nonlinear dispersive equations (local and global existence,
regularity, scattering theory) is vast and has been studied extensively by
many authors. Almost exclusively, the techniques developed so far restrict
to Cauchy problems with initial data in a Sobolev space, mainly because of the
crucial role played by the Fourier transform in the analysis of partial differential
operators. For a sample of results and a nice introduction to the field, we refer the
reader to Tao's monograph \cite{tao} and the references therein.

In this note, we focus on the Cauchy problem for the nonlinear Schr\"odinger equation (NLS),
the nonlinear wave equation (NLW), and the nonlinear Klein-Gordon equation (NLKG) in the realm of
modulation spaces. Generally speaking, a Cauchy data in a modulation space is rougher than any
given one in a fractional Bessel potential space and this low-regularity is desirable in many situations.
Modulation spaces were introduced by Feichtinger in the 80s \cite{Fei83} and have asserted
themselves lately as the ``right'' spaces in time-frequency analysis. Furthermore, they provide an excellent
substitute in estimates that are known to fail on Lebesgue spaces. This is not entirely surprising,
if we consider their analogy with Besov spaces, since modulation spaces arise essentially replacing dilation
by modulation.

The equations that we will investigate are:
\begin{equation}\label{nlschro}
\hspace*{-4cm}(NLS)\,\,\, i\frac{\partial{u}}{\partial{t}} + \Delta_{x}u + f(u)=0, \,  u(x, 0) = u_{0}(x),\,
\end{equation}
\begin{equation}\label{nlwave}
\hspace*{-.75cm}(NLW)\,\,\, \frac{\partial^{2}{u}}{\partial{t^{2}}} -\Delta_{x}u +
f(u)=0, \, u(x, 0) = u_{0}(x), \, \frac{\partial{u}}{\partial{t}}(x,
0)=u_{1}(x), \,
\end{equation}
\begin{equation}\label{nlklg}
(NLKG)\,\,\, \frac{\partial^{2}{u}}{\partial{t^{2}}} + (I-\Delta_{x})u
+ f(u)=0, \, u(x, 0)=u_{0}(x), \, \frac{\partial{u}}{\partial{t}}(x,
0) = u_{1}(x),
\end{equation}
where $u(x, t)$ is a complex valued function on $\rd\times \R$, $f(u)$ (the nonlinearity) is some scalar
function of $u$, and $u_0, u_1$ are complex valued functions on $\rd$.

The nonlinearities considered in this paper will be either power-like
\begin{equation}\label{powernl}
p_k(u)=\lambda |u|^{2k}u, k\in \N, \lambda\in\R,
\end{equation}
or exponential-like
 \begin{equation}\label{exponl}
 e_\rho (u)=\lambda (e^{\rho |u|^2}-1)u, \lambda, \rho\in\R.
 \end{equation}
 Both nonlinearities considered have the
advantage of being smooth. The corresponding equations having power-like
nonlinearities $p_k$ are
sometimes referred to as algebraic nonlinear (Schr\"odinger, wave, Klein-Gordon) equations.
The sign of the coefficient $\lambda$ determines the defocusing, absent, or
focusing character of the nonlinearity, but, as we shall see, this character will play no role in
our analysis on modulation spaces.

The classical definition of (weighted) modulation spaces that will be used throughout this work
is based on the notion of short-time Fourier transform (STFT). For $z=(x, \omega)\in\R^{2d}$, we let
$M_\omega$ and $T_x$ denote the operators of modulation and translation, and $\pi (z)=M_\omega T_x$ the general
time-frequency shift. Then, the STFT of $f$ with respect to a window $g$ is
$$V_g f(z)=\langle f, \pi (z)g\rangle .$$
Modulation spaces provide an effective way to measure the time-frequency concentration of a distribution
through size and integrability conditions on its STFT. For $s, t \in \R$ and $1\leq p, q\leq\infty$, we define the
weighted modulation space $\M{p,q}_{t, s}(\R^{d})$ to be the Banach space of all tempered distributions $f$ such that, for
a nonzero  smooth rapidly decreasing function $g \in \S(\R^{d})$, we have
$$\nm{f}{\M{p,q}_{t,s}} =
\biggparen{\int_{\R^{d}}\biggparen{\int_{\R^{d}}|V_{g}f(x,
\omega)|^{p}\, <x>^{tp}\, dx}^{q/p}\, <\omega>^{qs}\, d\omega}^{1/q} <\infty.$$ Here, we use the notation
$$<x>=(1+|x|^{2})^{1/2}.$$
This definition is independent of the choice of the window, in the
sense that different window functions yield equivalent
modulation-space norms. When both $s=t=0$, we will simply write
$\M{p, q}=\M{p, q}_{0, 0}$. It is well-known that the dual of a
modulation space is also a modulation space, $(\M{p, q}_{s,
t})'=\M{p', q'}_{-s, -t}$, where $p', q'$ denote the dual exponents
of $p$ and $q$, respectively. The definition above can be
appropriately extended to exponents $0<p, q\leq \infty$ as in the
works of Kobayashi \cite{ko1}, \cite{ko2}. More specifically, let
$\beta >0$ and $\chi \in \cS$ such that $\textrm{supp} \hat{\chi}
\subset \{|\xi|\leq 1\}$ and $\sum_{k\in \Z^{d}} \hat{\chi}(\xi -
\beta k) = 1, \, \forall \xi \in \R^d.$ For $0 <p,q\leq \infty$ and
$s>0$, the modulation space $\M{p,q}_{0,s}$ is the set of all
tempered distributions $f$ such that

\begin{equation}\label{komod}
\biggparen{\sum_{k \in \Z^{d}} \biggparen{\int_{\R^{d}}
|f \ast (M_{\beta k}\chi)(x)|^{p}\, dx}^{\tfrac{q}{p}} <\beta k>^{sq}}^{\tfrac{1}{q}} <\infty.
\end{equation}

When, $1\leq p, q \leq \infty$ this is an equivalent norm on
$\M{p,q}_{0,s}$, but when $0<p,q<1$ this is just a quasi-norm. We
refer to \cite{ko1} for more details. For another definition of the
modulation spaces for all $0<p,q\leq \infty$ we refer to
\cite{galsam, tri83}. For a discussion of the cases when $p$ and/or
$q=0$, see \cite{bgho}. These extensions of modulation spaces have
recently been rediscovered and many of their known properties
reproved via different methods by Baoxiang et all \cite{wang06},
\cite{bahud06}. There exist several embedding results between
Lebesgue, Sobolev, or Besov spaces and modulation spaces, see for
example \cite{Ok04}, \cite{toft04a}; also \cite{wang06},
\cite{bahud06}. We note, in particular, that the Sobolev space
$H^2_s$ coincides with $\M{2,2}_{0, s}$. For further properties and
uses of modulation spaces, the interested reader is referred to
Gr\"ochenig's book \cite{Gr01}.

The goal of this note is two fold: {\it to improve} some recent results of
Baoxiang, Lifeng and Boling \cite{wang06} on the local well-posedness of nonlinear equations
stated above, by allowing the Cauchy data to lie in any modulation space
$\M{p,1}_{0,s}$, $p>\tfrac{d}{d+1}$, $s\geq 0$, and {\it to simplify} the methods of proof by employing well-established tools from
time-frequency analysis. Ideally, one would like to adapt these methods to deal with global well-posedness
as well. We plan to address these issues in a future work.

In what follows, we assume that $d\geq 1, k\in \N, \tfrac{d}{d+1}< p\leq \infty$,
$\lambda, \rho\in \R$ and $s\geq 0$ are given. With $p_k$ and
$e_\rho$ defined by~\eqref{powernl} and~\eqref{exponl} respectively,
our main results are the following.

\begin{theorem}\label{lwpnls}
Assume that $u_0\in \M{p,1}_{0, s}(\rd)$ and $f\in \{p_k, e_\rho\}$. Then, there
exists $T^{*}=T^{*}(\nm{u_{0}}{\M{p,1}_{0, s}})$ such that~\eqref{nlschro}
has a unique solution $u \in C([0, T^{*}], \M{p,1}_{0, s}(\rd))$. Moreover,
if $T^{*}<\infty$, then $\displaystyle\limsup_{t\to T^{*}} \nm{u(\cdot,
t)}{\M{p,1}_{0, s}} = \infty.$
\end{theorem}

\begin{theorem}\label{lwpnw}
Assume that $u_0, u_1\in \M{p,1}_{0, s}(\rd)$ and $f\in \{p_k, e_\rho\}$. Then, there
exists $T^{*}=T^{*}(\nm{u_{0}}{\M{p,1}_{0,s}},\nm{u_{1}}{\M{p,1}_{0,s}})$ such that~\eqref{nlwave}
has a unique solution $u \in C([0, T^{*}], \M{p,1}_{0,s}(\rd))$. Moreover,
if $T^{*}<\infty$, then $\displaystyle\limsup_{t\to T^{*}} \nm{u(\cdot,
t)}{\M{p,1}_{0,s}} = \infty.$
\end{theorem}

\begin{theorem}\label{lwpnklg}
Assume that $u_0, u_1\in \M{p,1}_{0,s}(\rd)$ and $f\in \{p_k, e_\rho\}$. Then, there
exists $T^{*}=T^{*}(\nm{u_{0}}{\M{p,1}_{0, s}},\nm{u_{1}}{\M{p,1}_{0, s}})$ such that~\eqref{nlklg}
has a unique solution $u \in C([0, T^{*}], \M{p,1}_{0, s}(\rd))$. Moreover,
if $T^{*}<\infty$, then $\displaystyle\limsup_{t\to T^{*}} \nm{u(\cdot,
t)}{\M{p,1}_{0,s}} = \infty.$
\end{theorem}

\smallskip

\noindent{\bf Remark 1.} In Theorem~\ref{lwpnls} we can replace the (NLS) equation with
the following more general (NLS) type equation:
\begin{equation}
\hspace*{-4cm}(NLS)_\alpha \,\,\, i\frac{\partial{u}}{\partial{t}} + \Delta_{x}^{\alpha/2}u + f(u)=0, \,
u(x, 0) = u_{0}(x),\,
\end{equation}
for any $\alpha \in [0, 2]$ and $p\geq 1$. The operator $\Delta_x^{\alpha/2}$ is
interpreted as a Fourier multiplier operator (with $t$ fixed),
$\widehat \Delta_x^{\alpha/2}u (\xi, t)=|\xi|^\alpha \widehat u(\xi,
t)$. This strengthening will become evident from the preliminary
Lemma 1 of the next section.

\smallskip

\noindent{\bf Remark 2.} Theorems 1.1 and 1.2 of \cite{wang06} are
particular cases of Theorem~\ref{lwpnls} with $p=2$ and $s=0$.

\section{Fourier multipliers and multilinear estimates}

The generic scheme in the local existence theory is to establish
linear and nonlinear estimates on appropriate spaces that contain
the solution $u$. As indicated by the main theorems above, the
spaces we consider here are $\M{p,1}_{0,s}$, and we present the
appropriate estimates in the lemmas below. In fact, we will need
estimates on Fourier multipliers on modulation spaces. As proved in
\cite{bgor} and \cite{FeiNa}, a function $\sigma(\xi)$ is a symbol
of a bounded Fourier multiplier on $\M{p,q}$ for $1\leq p, q\leq
\infty$ if  $\sigma \in \wam$ (see the proofs of the following two
lemmas for a definition of this space). As we shall indicate below,
this condition can be naturally extended to give a sufficient
criterion for the boundedness of the Fourier multiplier operator on
$\M{p,q}_{0,s}$ for $0<p,q\leq \infty$ and $s\geq 0$. The notation
$A\lesssim B$ stands for $A\leq cB$ for some positive constant $c$
independent of $A$ and $B$.

\begin{lemma}\label{fmwlp}
Let $\sigma$ be a function defined on $\R^d$ and consider the Fourier multiplier operator $H_{\sigma}$ defined by $$H_{\sigma}f(x) =\int_{\R^{d}}\sigma(\xi)\, \hat{f}(\xi)\, e^{2\pi \xi \cdot x}\, d\xi.$$ Let $\chi \in \cS$ such that $\textrm{supp}\, \hat{\chi} \subset \{|\xi|\leq 1\}$. Let $d\geq 1$, $s\geq 0$, $0<q\leq \infty$, and $0<p<1$. If $\sigma \in \wamp(\R^{d})$, i.e., $$\nm{\sigma}{\wamp}= \sup_{n \in \Z^{d}}\nm{\sigma \cdot T_{\beta n} \chi}{\mathcal{F}L^{p}} <\infty$$ for $\beta>0$, then $H_{\sigma}$ extends to a bounded operator on $\M{p,q}_{0,s}(\R^{d})$.
\end{lemma}

\begin{proof}
We use the definition of the modulation spaces given by~\eqref{komod} (see also \cite{ko1}). In
particular, let $ \chi \in \cS$ such that $\textrm{supp}\,
\hat{\chi} \subset \{ |\xi|\leq 1\}$, and define  $g\in \cS$ by
$\hat{g}=\hat{\chi}^{2}.$ Denote $\tilde{g}(x)= \overline{g(-x)}.$
For $f\in \cS$,  $\beta
>0$, $k\in \Z^{d}$ and $x\in \R^d$ we have:
\begin{align*}
|H_{\sigma}f \ast(M_{\beta k} \tilde{g})(x)|& = |V_{g}H_{\sigma}f(x, \beta k)|\\
&= |\ip{\sigma \hat{f}}{M_{-x}T_{\beta k}\hat{g}} |\\
&= |\ip{\sigma \hat{f}}{M_{-x}T_{\beta k}\hat{\chi}^{2}} |\\
&\leq |\mathcal{F}^{-1}(\sigma \cdot T_{\beta k}\overline{\hat{\chi}})|\ast |\mathcal{F}^{-1}(\hat{f}\cdot T_{\beta k}\overline{\hat{\chi}})|(x)\\
& \leq |\mathcal{F}^{-1}(\sigma \cdot T_{\beta k}\overline{\hat{\chi}})|\ast |f\ast(M_{\beta k}\tilde{\chi})|(x).
\end{align*}
Now, observe that $\textrm{supp}\bigparen{\sigma \cdot T_{\beta
k}\overline{\hat{\chi}}} \subset \Gamma_k:=\beta k + \{|\xi|\leq
1\}$ and $\textrm{supp} \bigparen{ \hat{f}\cdot T_{\beta
k}\overline{\hat{\chi}}} \subset \Gamma_k.$ Moreover, by assumption  we know that $\sigma \in \wamp$ and so
$\mathcal{F}^{-1}(\sigma \cdot T_{\beta k}\overline{\hat{\chi}})
\in L^{p}$ and $f\ast(M_{\beta k}\tilde{\chi}) \in L^p$.
Consequently, by~\cite[Lemma 2.6]{ko1} we have the following estimate
$$\nm{H_{\sigma}f \ast(M_{\beta k} \tilde{g})}{L^{p}}
\leq C\, \nm{\mathcal{F}^{-1}(\sigma \cdot T_{\beta
k}\overline{\hat{\chi}})}{L^{p}} \nm{f\ast(M_{\beta
k}\tilde{\chi})}{L^{p}},$$ where $C$ is a positive constant that
depends only on the diameter of $\Gamma_k$ and $p$. Clearly, the
diameter of $\Gamma_k$ is independent of $k$, and this makes $C$ a
constant depending only on the dimension $d$ and the exponent $p$. Therefore, for $0<q \leq \infty$ we have
$$\nm{H_{\sigma}f}{\M{p,q}_{0,s}} \lesssim \sup_{k\in \Z^{d}} \nm{\mathcal{F}^{-1}(\sigma \cdot
T_{\beta k}\overline{\hat{\chi}})}{L^{p}}\, \nm{f}{\M{p,q}_{0,s}}= \nm{\sigma}{\wamp}\, \nm{f}{\M{p,q}_{0,s}}.$$
The result then follows from the density of $\cS$ in $\M{p,q}_{0,s}$ for $p, q <\infty$; see \cite[Theorem 3.10]{ko1}.
\end{proof}

\smallskip

We are now ready to state and prove the boundedness of Fourier
multipliers that will be needed in establishing our main results.

\begin{lemma}\label{weimult}
Let $d\geq 1$, $s\geq 0$,  and $0<q\leq \infty$ be
given. Define $m_{\alpha}(\xi)=e^{i|\xi|^{\alpha}}$.  
If $1\leq p \leq \infty$ and $\alpha \in [0, 2]$,  then the Fourier multiplier operator $H_{m_{\alpha}}$ extends to a bounded operator on $\M{p,q}_{0, s}(\R^{d})$.

Moreover, If $\alpha \in \{1, 2\}$ and $\tfrac{d}{d+1} <p \leq \infty$, then the Fourier multiplier operator $H_{m_{\alpha}}$ extends to a bounded operator on $\M{p,q}_{0, s}(\R^{d})$.
\end{lemma}

\begin{proof}
First, we prove the result when $1\leq p \leq \infty$, and $0<q\leq
\infty$. Let $g \in \S(\R^{d})$ and define $\chi \in \S$ by
$\hat{\chi}=g^2$. For $f \in \S$, we have
\begin{align}\nonumber
|V_{\chi}&H_{m_{\alpha}}f (x, \xi)|\\
&= \bigg{|} \int_{\R^{d}} m_{\alpha}(t)\, \hat{f}(t)\, e^{2\pi i x\cdot t}\,
\overline{\hat{\chi}(t - \xi)}\, dt\bigg{|}\notag \\
& = \tfrac{1}{<\xi>^{s}}\bigg{|} \int_{\R^{d}} m_{\alpha}(t)\,
\overline{T_{\xi}g(t)}\,  <t>^{s}\, \hat{f}(t)\, \tfrac{<\xi>^{s}}{<t>^{s}\, <t-\xi>^{N}}\,  <t-\xi>^{N}\,
\overline{g(t-\xi)}\, e^{2\pi i x\cdot t}\, dt \bigg{|} \notag \\
&= \tfrac{1}{<\xi>^{s}}\bigg{|} \int_{\R^{d}} m_{\alpha}(t)\,
\overline{T_{\xi}g(t)}\,\phi_N(\xi, t)\, \widehat{<D>^{s}f}(t)\, T_{\xi}\overline{g_{N}(t)}\, dt\bigg{|}\notag \\
& =\tfrac{1}{<\xi>^{s}}\bigg{|}\mathcal{F} \biggparen{m_{\alpha}\cdot T_{\xi}g\, \phi_N(\xi, \cdot)\,
\widehat{<D>^{s}f}\, T_{\xi}\overline{g_{N}}}(-x)\bigg{|}\notag\\
& =\tfrac{1}{<\xi>^{s}}\bigg{|} \mathcal{F}(m_{\alpha}\cdot
T_{\xi}g)\ast \mathcal{F}_{2}(\phi_N(\xi, \cdot))\ast \mathcal{F}(
\widehat{<D>^{s}f}\cdot T_{\xi}\overline{g_{N}}) (-x)\bigg{|}\notag,
\end{align}
where $N>0$ is an integer to be chosen later, $g_{N}(t) =
<t>^{N}\overline{g}(t)$, $\phi_N (\xi, t) = \tfrac{<\xi>^{s}}{<t>^{s}\,
<t-\xi>^{N}}$, and $<D>^s$ is the Fourier
multiplier defined by $\widehat{<D>^{s}f}(\xi) =
<\xi>^{s}\hat{f}(\xi)$.  We also denote by $\Phi_{2, N}(\xi, \cdot):=\mathcal{F}_{2}(\phi_N (\xi, \cdot))$
the Fourier transform in the second variable of $\phi_N (\xi, \cdot)$

We can therefore estimate the weighted modulation norm of $H_{m_\alpha}f$ as follows:
\begin{align}\nonumber
&\nm{H_{m_{\alpha}}f}{\M{p,q}_{0, s}}\\
&= \biggparen{\int_{\R^{d}}\biggparen{\int_{\R^{d}}|V_{\chi}f(x, \xi)|^{p}\, dx}^{q/p}\,
<\xi>^{qs}}^{1/q} \notag \\
&= \biggparen{\int_{\R^{d}} \biggparen{\int_{\R^{d}} \bigg{|}
\mathcal{F}(m_{\alpha}\cdot T_{\xi}g)\ast \Phi_{2, N}(\xi, \cdot)\ast \mathcal{F}
( \widehat{<D>^{s}f}\cdot T_{\xi}\overline{g_{N}}) (-x)\bigg{|}^{p}dx}^{q/p}d\xi}^{1/q}\notag \\
& \leq \biggparen{\int_{\R^{d}}\nm{\mathcal{F}^{-1}(m_{\alpha}\cdot T_{\xi}g)}{L^{1}}^{q}\,
\nm{\Phi_{2, N}(\xi, \cdot)}{L^{1}}^{q}\, \nm{\mathcal{F}( \widehat{<D>^{s}f}
\cdot T_{\xi}\overline{g_{N}})}{L^{p}}^{q}\, d\xi}^{1/q}\notag \\
& \leq \sup_{\xi \in \R^{d}} \nm{\mathcal{F}^{-1}(m_{\alpha}\cdot T_{\xi}g)}{L^{1}}\,
\sup_{\xi \in \R^{d}} \nm{\Phi_{2, N}(\xi, \cdot)}{L^{1}}\,
\biggparen{\int_{\R^{d}} \nm{\mathcal{F}^{-1}( \widehat{<D>^{s}f}\cdot T_{\xi}
\overline{g_{N}})}{L^{p}}^{q}\, d\xi}^{1/q}\notag \\
& \leq \sup_{\xi \in \R^{d}} \nm{\mathcal{F}(m_{\alpha}\cdot
T_{\xi}g)}{L^{1}}\, \sup_{\xi \in \R^{d}} \nm{\Phi_{2, N}(\xi,
\cdot)}{L^{1}}\, \nm{f}{\M{p,q}_{0, s}}. \label{multconv2}
\end{align}
Now, it follows from \cite[Lemma 8]{bgor} that, for $\alpha \in [0,2]$,
$$ \sup_{\xi \in \R^{d}} \nm{\mathcal{F}^{-1}(m_{\alpha}\cdot T_{\xi}g)}{L^{1}}:=
\nm{m_{\alpha}}{\wam} <\infty.$$
Moreover (see, for example, \cite[Lemma 3.1]{toft04a} or \cite[Lemma
2.1]{toft04b}), we can select a sufficiently large $N>0$ such that
$$\sup_{\xi \in \R^{d}} \nm{\Phi_{2, N}(\xi, \cdot)}{L^{1}} \leq
\int_{\R^{d}} \sup_{\xi \in \R^{d}} |\Phi_{2, N}(\xi, x))| dx <
\infty.$$
Hence, using ~\eqref{multconv2}, we get
$$\nm{H_{m_{\alpha}}f}{\M{p,q}_{0, s}}\leq C_{\alpha}\nm{f}{\M{p,q}_{0, s}}.$$

To prove the second part of the result we shall use
Lemma~\ref{fmwlp}. In particular, we need to show that for $\alpha
\in \{1,2\}$ and $\tfrac{d}{d+1}<p<1$, $m_{\alpha} \in \wamp$. This,
however, follows by straightforward adaptations of the proofs of
\cite[Theorems 9 and 11]{bgor}, which we leave to the interested
reader.
\end{proof}

\smallskip
In analogy to the proof of the previous lemma, we can prove the following weighted version of
\cite[Theorem 16]{bgor}.

\begin{lemma}\label{weimult2}
Let $d\geq 1$, $s\geq 0$, $\tfrac{d}{d+1}< p\leq \infty$ and $0< q\leq\infty$ be given, and let $m^{(1)}(\xi) = \tfrac{\sin (|\xi|)}{|\xi|}$
and $m^{(2)}(\xi)= \cos (|\xi|),$ for $\xi \in \R^{d}$. Then, the Fourier multiplier operators $H_{m^{(1)}}, H_{m^{(2)}}$ can be extended as
bounded operators on $\M{p,q}_{0, s}.$
\end{lemma}

A ``smooth'' version of Lemma \ref{weimult2} is obtained by replacing $|\xi|$ with $<\xi>$.

\begin{lemma}\label{multklg}
Let $d\geq 1$, $s\geq 0$, $\tfrac{d}{d+1} < p\leq \infty$ and $0 < q\leq\infty$ be given, and let $m(\xi) =e^{i<\xi>}$,
$m^{(1)}(\xi) = \tfrac{\sin (<\xi>)}{<\xi>}$ and $m^{(2)}(\xi)= \cos (<\xi>),$ for $\xi \in \R^{d}$.
Then, the Fourier multiplier operators
$H_m, H_{m^{(1)}}, H_{m^{(2)}}$ can be extended as bounded operators on $\M{p,q}_{0, s}.$
\end{lemma}

\begin{proof}
It is clear that $m, m^{(1)}, m^{(2)}$ are $\mathcal{C}^{\infty}(\R^{d})$ functions and that all their
derivatives are bounded. Therefore, $m, m^{(1)}, m^{(2)} \in \mathcal{C}^{d+1}(\R^{d}) \subset \M{\infty, 1}(\R^{d})\subset \wam(\R^{d})$ \cite{Gr01, Ok04}. Thus, for $1\leq p \leq \infty$, and $0<q\leq \infty$ the result follows from \cite{bgor} and Lemma~\ref{weimult}. For $\tfrac{d}{d+1}<p<1$ and $0<q\leq \infty$, it can be showed that $m, m^{(1)}, m^{(2)}\in \mathcal{C}^{d+1}(\R^{d}) \subset \wamp(\R^{d})$. Indeed, this follows from obvious modifications to the proof of the embedding  $\mathcal{C}^{d+1}(\R^{d})
\subset \M{\infty, 1}(\R^{d}) \subset \wam(\R^{d})$ \cite{Gr01, Ok04}.  Furthermore, if we modify, for example, the multiplier $m$ to
$m_t(\xi)=e^{it<\xi>}$, $t\in\R$, we have for $\tfrac{d}{d+1} < p\leq 1$
\begin{equation}\label{multiplierklg}
\nm{m_t}{\wamp} \leq (1 + |t|)^{d+1},
\end{equation}
and similar estimates hold for modified multipliers $m^{(1)}_t$ and $m^{(2)}_t$.
\end{proof}

Finally, we state a crucial multilinear estimate that will be used in our proofs. Although the
estimate will be needed only in the particular case of a product of functions (see Corollary 1), we present it
here in its full generality that applies to multilinear pseudodifferential operators.

An $m$-linear pseudodifferential operator is defined \`a priori through its
(distributional) symbol $\sigma$ to be the mapping $T_\sigma$ from the
$m$-fold product of Schwartz spaces $\S \times \cdots \times \S$
into the space $\di$ of tempered distributions given by
the formula
\begin{align} \label{mpdo}
T_\sigma(u_1, & \dots, u_m)(x) \notag \\
& = \int_{\rd{m}} \sigma(x,\xi_1, \dots, \xi_m) \,
  \hat{u_1}(\xi_1) \, \cdots \, \hat{u}_m(\xi_m) \,
  e^{2 \pi i x \cdot (\xi_1 + \cdots + \xi_m)} \, d\xi_1 \, \cdots \, d\xi_m,
\end{align}
for $u_1, \dots, u_m \in \S$. The pointwise product $u_1 \cdots u_m$ corresponds to the case
$\sigma =1$.

\begin{lemma}\label{multimod}
If $\sigma \in \M{\infty,1}_{0, s}(\R^{(m+1)d})$,
then the $m$-linear pseudodifferential operator $T_\sigma$
defined by~\eqref{mpdo} extends to a bounded operator from
$\M{p_1,q_1}_{0,s} \times \cdots \times \M{p_m,q_m}_{0,s}$ into $\M{p_0, q_0}_{0, s}$
when
$\frac{1}{p_1} + \cdots + \frac{1}{p_m} = \frac{1}{p_0}$,
$\frac{1}{q_1} + \cdots + \frac{1}{q_m} = m-1 + \frac{1}{q_0}$,
and $0<p_i\leq\infty, 1\leq q_i \leq \infty$ for $0 \leq i \leq m$.
\end{lemma}

This result is a slight modification of \cite[Theorem 3.1]{bgho}. Its proof proceeds
along the same lines, and therefore it is omitted here. Note that if $\sigma\in\M{\infty, 1}_{0, s}$, and
we pick $u_1=\dots =u_m=u$ (some of them could be equal to $\bar u$ since the modulation norm is preserved),
$p_1=\dots =p_m=mp,$  $0< p\leq \infty$, and $q_1=\dots=q_m=1$ we have
\begin{equation}\label{multi2}
\nm{T_\sigma (u, \dots, u)}{\M{p, 1}_{0, s}}\lesssim \nm{u}{\M{mp, 1}_{0,s}}^m\lesssim \nm{u}{\M{p, 1}_{0,s}}^m,
\end{equation}
where we used the obvious embedding $\M{p, 1}_{0,s}\subseteq \M{mp, 1}_{0,s}.$ The notation $A\lesssim B$ stands for
$A\leq cB$ for some positive constant $c$ independent of $A$ and $B$. In particular, if we select
$\sigma=1$ (the constant function 1), then $\sigma\in \M{\infty, 1}_{0,s} \subset \M{\infty,1}$, and we obtain

\begin{corollary}\label{prod1}
Let $0<p\leq \infty$. If $u\in \M{p, 1}_{0,s}$, then $u^m\in \M{p, 1}_{0,s}$. Furthermore,
$$\nm{u^m}{\M{p, 1}_{0,s}}\lesssim \nm{u}{\M{p, 1}_{0,s}}^m.$$
\end{corollary}
This is of course just a particular case of the more general multilinear estimate
\begin{equation}\label{prod2}
\bignm{\displaystyle\prod_{i=1}^m u_i}{\M{p_0, q_0}_{0,s}}\lesssim \displaystyle\prod_{i=1}^m\nm{u_i}{\M{p_i, q_i}_{0,s}},
\end{equation}
where the exponents satisfy the same relations as in Lemma 1. When we consider the
power nonlinearity $f(u)=p_k(u)=\lambda |u|^{2k}u=\lambda u^{k+1}\bar u^{k}$, Corollary \ref{prod1} becomes

\begin{corollary}\label{prod3}
Let $0<p\leq \infty$. If $u\in \M{p, 1}_{0,s}$, then $p_k(u)\in \M{p, 1}_{0,s}$. Furthermore,
$$\nm{p_k(u)}{\M{p, 1}_{0,s}}\lesssim \nm{u}{\M{p, 1}_{0,s}}^{2k+1}.$$
\end{corollary}
For a different proof of the estimate in Corollary \ref{prod3},
see \cite[Corollary 4.2]{wang06}. It is important to note that the previous estimate allows us to control
the exponential nonlinearity $e_\rho$ as well. Indeed, since
$$e_\rho (u)=\lambda (e^{\rho |u|^{2}} -1) u =
\sum_{k=1}^{\infty} \tfrac{\rho^{k}}{k!}p_k (u),$$
if we now apply the modulation norm on both sides and use the triangle inequality, we arrive at

\begin{corollary}\label{prod4}
Let $0<p\leq \infty$. If $u\in \M{p, 1}_{0,s}$, then $e_\rho (u)\in \M{p, 1}_{0,s}$. Furthermore,
$$\nm{e_\rho (u)}{\M{p, 1}_{0,s}}\lesssim \nm{u}{\M{p, 1}_{0,s}}(e^{|\rho| \nm{u}{\M{p, 1}_{0,s}}^2}-1).$$
\end{corollary}

\section{Proofs of the main results}

We are now ready to proceed with the proofs of our main theorems. We will only prove our results
for the power nonlinearities $f=p_k$, by making use of Corollary \ref{prod3}. The case of exponential
nonlinearity $f=e_\rho$ is treated similarly, by now employing Corollary \ref{prod4}. In all that follows we assume
that  $u: [0, T)\times \rd \to \C$ where $0< T \leq \infty$ and that $f(u)=p_k(u)=\lambda |u|^{2k}u.$

\subsection{The nonlinear Schr\"odinger equation: Proof of Theorem \ref{lwpnls}}
We start by noting that ~\eqref{nlschro} can be written in the equivalent form
\begin{equation}\label{eqschro}
u(\cdot , t) = S(t)u_{0}  - i \mathcal{A}f(u)
\end{equation} where
\begin{equation}\label{nlssemi}
S(t) = e^{i t \Delta}, \, \, \, \, \, \mathcal{A} =
\int_{0}^{t}S(t-\tau) \cdot \, d\tau.
\end{equation}

Consider now the mapping
$$\mathcal{J} u = S(t)u_{0} - i\int_{0}^{t}S(t-\tau) (p_k(u))(\tau)\, d\tau.$$ It follows from
Lemma \ref{weimult} (see also \cite[Corollary 18]{bgor}) that
\begin{equation}\nonumber
\nm{S(t)u_{0}}{\M{p,1}_{0,s}} \leq C\, (t^2 + 4\pi^2)^{d/4}\,\nm{u_{0}}{\M{p,1}_{0,s}},
\end{equation}
where  $C$ is a universal constant depending only on $d$. Therefore,
\begin{equation}\label{step1}
\nm{S(t)u_{0}}{\M{p,1}_{0,s}} \leq C_T\, \nm{u_{0}}{\M{p,1}_{0,s}},
\end{equation}
where $C_T = \displaystyle\sup_{t \in [0, T)}C\, (t^2 + 4\pi^2)^{d/4}.$  Moreover, we have
\begin{align}
\bignm{\int_{0}^{t}S(t-\tau) (p_k(u))(\tau)\, d\tau}{\M{p,1}_{0,s}} & \leq \, \int_{0}^{t} \nm{S(t-\tau)
(p_k(u))(\tau)}{\M{p,1}_{0,s}}\, d\tau \notag \\
& \leq T \,C_T \, \, \sup_{t \in [0, T]} \nm{p_k(u)(t)}{\M{p,1}_{0,s}}. \label{est2}
\end{align}

By using now Corollary \ref{prod3}, we can further estimate in ~\eqref{est2} to get
\begin{equation}\label{est3}
\bignm{\int_{0}^{t}S(t-\tau) (p_k(u))(\tau)\,
d\tau}{\M{p,1}_{0,s}}  \lesssim C_{T}\, T\,
\nm{u(t)}{\M{p, 1}_{0,s}}^{2k+1}.
\end{equation}
Consequently, using~\eqref{step1} and~\eqref{est3} we have
\begin{equation}\label{est4}
\nm{\mathcal{J}u }{C([0,T], \M{p,1}_{0,s})} \leq C_{T} (\nm{u_{0}}{\M{p,1}_{0,s}} + cT \, \nm{u}{\M{p,1}_{0,s}}^{2k+1}),
\end{equation}
for some universal positive constant $c$. We are now in the position of using a standard contraction argument
to arrive to our result. For completeness, we sketch it here. Let $\mathbf{B}_M$ denote the closed ball of radius $M$
centered at the origin in the space $C([0, T], \M{p, 1}_{0,s})$. We claim that
$$\mathcal{J}:\mathbf{B}_M\rightarrow \mathbf{B}_M,$$
for a carefully chosen $M$. Indeed, if we let $M=2C_T\nm{u_0}{\M{p, 1}_{0,s}}$ and $u\in\mathbf{B}_M$, from (\ref{est4}) we
obtain
$$\nm{\mathcal{J}u }{C([0,T], \M{p,1}_{0,s})} \leq \frac{M}{2}+cC_TTM^{2k+1}.$$
Now let $T$ be such that $cC_TTM^{2k}\leq 1/2$, that is, $T\leq \tilde T(\nm{u_0}{\M{p, 1}_{0,s}})$. We obtain
$$\nm{\mathcal{J}u }{C([0,T], \M{p,1}_{0,s})} \leq \frac{M}{2}+\frac{M}{2}=M,$$
that is $\mathcal{J}u\in\mathbf{B}_M$. Furthermore, a similar argument gives
$$\nm{\mathcal{J}u-\mathcal{J}v}{C([0,T], \M{p,1}_{0,s})}\leq \frac{1}{2}\nm{u-v}{C([0,T], \M{p,1}_{0,s})}.$$ This last
estimate follows in particular from the following fact:
$$p_{k}(u)(\tau) -p_{k}(v)(\tau) = \lambda (u-v)|u|^{2k}(\tau) +
\lambda v (|u|^{2k} -|v|^{2k})(\tau).$$ Therefore, using Banach's
contraction mapping principle, we conclude that $\mathcal{J}$ has a
fixed point in $\mathbf{B}_M$ which is a solution of
(\ref{eqschro}); this solution can be now extended up to a maximal
time $T^{*}(\nm{u_0}{\M{p, 1}_{0,s}})$. The proof is complete.

\subsection{The nonlinear wave equation: Proof of Theorem \ref{lwpnw}}

Equation ~\eqref{nlwave} can be written in the equivalent form

\begin{equation}\label{eqnlwave}
u(\cdot, t) = \tilde K(t)u_{0} + K(t) u_{1} -\mathcal{B}f(u)
\end{equation} where
\begin{equation}\label{nlwavesemi}
K(t) =\tfrac{\sin(t \sqrt{ -\Delta})}{\sqrt{ - \Delta}}, \, \, \tilde K(t)
= \cos (t \sqrt{ -\Delta}), \,\, \mathcal{B}=\int_{0}^{t}K(t -\tau)\cdot
d\tau
\end{equation}

Consider the mapping
$$\mathcal{J} u = \tilde K(t) u_{0} + K(t) u_{1} -\mathcal{B}f(u).$$
Recall that $f=p_k$. If we now use Lemma \ref{weimult2} (see also \cite[Corollary
21]{bgor}) for the first two inequalities below and Corollary \ref{prod3} for the last estimate, we can write

\begin{equation}\label{est5}
\left\{\begin{array}{r@{,}l}
\nm{\tilde K(t)u_{0}}{\M{p,1}_{0,s}} \leq C_{T}\nm{u_0}{\M{p,1}_{0,s}}& \\
\nm{K(t)u_{1}}{\M{p,1}_{0,s}} \leq C_{T}\nm{u_1}{\M{p,1}_{0,s}}&\\
\nm{\mathcal{B}f(u)}{\M{p,1}_{0,s}}\leq cT\,
C_{T} \nm{u}{\M{p,1}_{0,s}}^{2k+1}&
\end{array}\right.
\end{equation} where $c$ is some universal positive constant. The constants $T$ and $C_T$ have the same meaning as before.
The standard contraction mapping argument applied to $\mathcal{J}$ completes the proof.

\subsection{The nonlinear Klein-Gordon equation: Proof of Theorem \ref{lwpnklg}}
The equivalent form of  equation ~\eqref{nlklg} is

\begin{equation}\label{eqnlklg}
u(\cdot, t) = \tilde K(t) u_{0} + K(t) u_{1} + \mathcal{C}f(u)
\end{equation} where now
\begin{equation}\label{nlklgsemi}
K(t) = \tfrac{\sin t (I-\Delta)^{1/2}}{(I-\Delta)^{1/2}}, \, \,
\tilde K(t) = \cos t (I-\Delta)^{1/2} ,\, \, \mathcal{C} = \int_{0}^{t}
K(t-\tau)\cdot d\tau.
\end{equation}

Consider the mapping $$\mathcal{J} u = \tilde K(t) u_{0} + K(t) u_{1} +
\mathcal{C}f(u).$$ Using Lemma \ref{multklg} and the notations above, we can write

\begin{equation}\label{est6}
\left\{\begin{array}{r@{,}l}
\nm{\tilde K(t)u_{0}}{\M{p,1}_{0,s}} \leq C_{T}\nm{u_0}{\M{p,1}_{0,s}}& \\
\nm{K(t)u_{1}}{\M{p,1}_{0,s}} \leq C_{T}\nm{u_1}{\M{p,1}_{0,s}}& \\
\nm{\mathcal{C}f(u)}{\M{p,1}_{0,s}}\leq cT\, C_{T}
\nm{u}{\M{p,1}_{0,s}}^{2k+1}&
\end{array}\right.
\end{equation}
The standard contraction mapping argument applied to $\mathcal{J}$ completes the proof.

\end{document}